\documentclass[12pt]{amsart}

\usepackage{amsfonts, amstext, amsmath, amsthm, amscd, amssymb, graphicx, epstopdf}
\usepackage{epsfig, graphics, psfrag, enumerate}
\usepackage{color}
\usepackage{verbatim}
\let\oldmarginpar\marginpar
\renewcommand\marginpar[1]{\oldmarginpar[\raggedleft\footnotesize #1]%
{\raggedright\footnotesize #1}}

\newtheorem{theorem}{Theorem}

\newtheorem{corollary}[theorem]{Corollary}

\newtheorem{conjecture}[theorem]{Conjecture}

\theoremstyle{definition}
\newtheorem{remark}[theorem]{Remark}

\newcommand{\NN}{{\mathbb{N}}}

\newcommand{\QQ}{{\mathbb{Q}}}

\newcommand{\bdy}{{\partial}}

\newcommand{\K}{{K_{p, q}}}

\newcommand{\abs}[1]{{\left\vert #1 \right\vert}}

\newcommand{\GA}{{\mathbb{G}_A}}

\newcommand\no[1]{}

\newtheorem*{namedtheorem}{\theoremname}
\newcommand{\theoremname}{testing}

\def\be { \begin{equation} }
\def\ee { \end{equation} }

\begin{document}

\title[]{The Strong Slope Conjecture and torus knots}

\author[]{Efstratia Kalfagianni}
\address{Department of Mathematics, Michigan State University, East Lansing, MI, 48824}
\email{kalfagia@math.msu.edu}
\bigskip

\bigskip

\begin{abstract} We show that the strong slope conjecture implies that the degree of the colored Jones polynomial detects all torus knots. As an application we obtain that an adequate knot that has the same colored Jones polynomial degrees as a torus knot must be a $(2,q)$-torus knot.
\bigskip
\bigskip

\bigskip

\bigskip
\end{abstract}

\bigskip

\bigskip
\thanks {\today}
\thanks{Supported in part by NSF grant DMS-1708249}

\maketitle

\section{Introduction} 
The colored Jones polynomial of a knot $K$ is a collection 
of Laurent polynomials 
$ \{ J_K(n):=J_{K}(n, t)\}_{n=1}^\infty $ 
in a variable $t$ such that we have $J_{K}(1, t)=1$ and $J_{K}(2,t)$ is the classical Jones polynomial.
We will use the normalization 
$$J_{\text{unknot}}(n) = \frac{t^{n/2}- t^{-n/2}}{t^{1/2} -t^{-1/2}}.$$ 

For a knot  $K \subset S^3$, let $d_+[J_{K}(n)]$ and  $d_-[J_{K}(n)]$ denote  the maximal and minimal degrees of $J_{K}(n)$ in $t$, respectively.
The strong slope conjecture \cite{ga-slope, Effie-Anh-slope} asserts that 
the degrees  $d_+[J_{K}(n)]$ and  $d_-[J_{K}(n)]$  contain information about essential surfaces in knot exteriors.
The purpose of this note is to show the following.
  
\begin{theorem} \label{main}Suppose that $K$ is  a knot that satisfies the strong slope conjecture and let $T_{p,q}$ denote the $(p,q)$-torus knot.
If $d_+[J_{K}(n)]=d_+[J_{T_{p,q}}(n)]$ and  $d_-[J_{K}(n)]=d_-[J_{T_{p,q}}(n)]$, for all $n$,  then, up to orientation change of the knot,  $K$ is isotopic to  $T_{p,q}$.
\end{theorem}

 The conjecture has been proved for several classes of knots and in particular for adequate knots \cite{FKP}.
 Using this we have the following.
 
 \begin{corollary} \label{adequate} Suppose that $K$ is an adequate knot such that  $d_+[J_{K}(n)]=d_+[J_{T_{p,q}}(n)]$ and  $d_-[J_{K}(n)]=d_-[J_{T_{p,q}}(n)]$, for all $n$.
 Then we have $|p|=2$.
 \end{corollary}
\vskip 0.07in

\section{Background}
It is known that 
the degrees $ d_+[J_{K}(n)] $ and $ d_-[J_{K}(n)] $ are quadratic { quasi-polynomials} \cite{ga-quasi}.
This means that, given a knot $K$, there is $n_K\in \NN$ such that  for all $n>n_K$ we have
 $$d_+[J_{K}(n)] =  a(n) n^2 + b(n) n  + c(n)\ \ \  {\rm and}  \ \ \   d_-[J_{K}(n)] =  a^{*}(n) n^2 + b^{*}(n) n  + c^{*}(n),$$
 where the coefficients are periodic functions from $\NN $ to $\QQ$ with integral period. By taking  the least common multiple of the periods
 of these coefficient  functions we get
 a common period. This common period of the coefficient functions is called the {\em Jones period} of $K$.

For a sequence $\{x_n\}$, let $\{x_n\}'$ denote the set of its cluster points. 
 The elements of the sets 
$$js_K:= \left\{ 4n^{-2}d_+[J_K(n)]  \right\}' \quad
 \mbox{and} \quad js^*_K:= \left\{ 4n^{-2}d_-[J_K(n)] \right\}' $$
 are called {\em Jones slopes} of $K$.

To continue, 
let  $\ell d_+[J_K(n)]$ and $\ell d_-[J_K(n)]$ denote the linear terms of $d_+[J_K(n)]$ and  $d_-[J_K(n)]$, respectively.
Also let 

$$jx_K:= \left\{ 2n^{-1}\ell d_+[J_K(n)]  \right\}'  \quad
 \mbox{and}\quad  jx^*_K:= \left\{ 2n^{-1}\ell d_-[J_K(n)]  \right\}'.$$

Given a knot $K\subset S^3$, let
  $n(K)$ denote a tubular neighborhood of
$K$ and let $M_K:=\overline{ S^3\setminus n(K)}$ denote the exterior of
$K$. Let $\langle \mu, \lambda \rangle$ be the canonical
meridian--longitude basis of $H_1 (\bdy n(K))$.   A properly embedded surface $(S, \bdy S) \subset (M_K,
\bdy n(K))$, is called  {\em essential }if it's $\pi_1$-injective and it is  not a boundary parallel annulus.
An element $\alpha/\beta \in
{\QQ}\cup \{ 1/0\}$, where $\alpha$ and $\beta$ are relatively prime integers,  is called a \emph{boundary slope} of $K$ if there
is an essential surface $(S, \bdy S) \subset (M_K,
\bdy n(K))$, such that each component of  $\bdy S$ represents $\alpha \mu + \beta \lambda \in
H_1 (\bdy n(K))$. The longitude $\lambda$ of every knot bounds an essential orientable surface in the exterior of $K$.  Thus
 $0=0/1$ is a boundary slope of every knot in $S^3$.
Hatcher showed that every knot $K \subset S^3$
has finitely many boundary slopes \cite{hatcher}. 

For a surface  $(S, \bdy S) \subset (M_K,
\bdy n(K))$ we will use the notation $\abs{\partial S}$ to denote the number of boundary components of $S$.

Garoufalidis conjectured \cite[Conjecture 1.2]{ga-slope}, that the Jones slopes of any knot $K$ are 
 boundary slopes.
The following statement, which  is a refinement of the original conjecture, was 
stated by the author and Tran in \cite[Conjecture 1.6]{Effie-Anh-slope}.

\begin{conjecture}  {\rm ({Strong slope conjecture})} \label{SSC}
\begin{itemize}  
\item Given a Jones slope $\alpha/\beta\in js_K$, with $\beta>0$ and $\gcd(\alpha, \beta)=1$, there is an essential surface $S$ in $M_K$ such that each component of $\partial S$ has slope $\alpha/\beta$  and
we have  $\displaystyle{  \frac{\chi(S)}{{\abs{\partial S} \beta}}} \in  jx_K$.

\item  Given  a Jones slope $\alpha^{*}/\beta^{*}\in  js^*_K$,  with $\beta^{*}>0$ and $\gcd(\alpha^{*}, \beta^{*})=1$,
there is an essential surface $S^{*}$ in $M_K$ such that each component of $\partial S^{*}$ has slope $\alpha^{*}/\beta^{*}$ and we have $\displaystyle{ - \frac{\chi(S^{*})}{{\abs{\partial S^{*}} \beta^{*}}}} \in  jx^*_K$.

\end{itemize}

\end{conjecture}

Conjecture \ref{SSC} has been proved for several classes of knots including Montesinos knots, adequate knots, graph knots and  cables.
The reader  is referred to \cite{ MoTa, Effie-Anh-slope, Takata, TakataMB, LeeVeen, LeeGarVeen} and references therein.

 \section {Proofs}
 
 We now prove Theorem \ref{main} and Corollary  \ref{adequate}.
  \smallskip
  
 {\it Proof of Theorem \ref{main}.} \
  First suppose that  $T_{p,q}$ is the trivial knot; that is $p=\pm 1$, or $q=\pm 1$. Then by the normalization of $J_{K}(n)$ fixed earlier, we have
  $4d_+[J_{K}(n)]=2(n-1)$. Since $K$ satisfies Conjecture \ref{SSC}  we have a surface  $S$ that is essential in exterior $M_K$, $\partial S$ has slope zero and we have
  $\displaystyle{{\chi(S)\over {{\abs{\partial S}} }}}=1$.  This implies that $S$ is a collection of discs and hence $K$ is the unknot.
 
 Assume now that $T_{p,q}$ is not the unknot. By  the classification of torus knots (see, for example, \cite[Theorem 3.29]{burde-zieschang:knots})
 we can assume that $\abs{p}<\abs{q}$ and $q>0$.

 For a knot  $K$ let $K^{*}$ denote the mirror image of $K$. It is known that
  $d_{\pm}[J_{K^{*}}(n)]=- d_{\mp}[J_{K}(n)]$, and that the colored Jones polynomial is insensitive to orientation change of knots.
  Let $\K$ denote the $(p, q)$-cable of $K$, where $p$ and $q$ are relatively prime integers.
We have $(K_{p,q})^{*}= K^{*}_{-p,q}$ and hence we have
 $$d_{\mp}[J_{K^{*}_{-p,q}}(n)]=- d_{\pm}[J_{K_{p,q}}(n)].$$
 In particular, $T_{p,q}$ is the $(p, q)$-cable of the unknot and  $T_{-p,q}$ is the mirror image 
 of $T_{p,q}$ while $T_{-p,-q}$ is $T_{p,q}$ with the orientation reversed.
 Without loss of generality we will assume that $p<0<q$.

  The degrees  $d_{\pm}[J_{T_{p ,q}}(n)]$ were calculated in  \cite{ga-slope}:
  we have
 $$4d_-[J_{K}(n)]=4d_-[J_{T_{p,q}}(n)]= pqn^2+d(n),$$
where $d(n)=-pq-\big( 1+(-1)^{n}\big)(p-2)(q-2)/2$, and we have
 $$4d_+[J_{K}(n)]=4d_+[J_{T_{p ,q}}(n)]= 2 (pq-p+q)n-2(pq-p+q).$$
 
 Note that the values of the quadratic and linear terms of the above formulae can also be obtained by the argument in the proof of \cite[Theorem 3.9]{Effie-Anh-slope}
 using the fact that $T_{p ,q}$ is  a cable of the unknot.
 
 Since the quadratic and linear coefficients of $ d_{\pm}[J_{K}(n)]$ are independent of $n$, each of the sets  $ js_K,   jx_K,  js^*_K,   jx^*_K$ consists of exactly one element (i.e. the corresponding coefficient).
 
Since $K$ is assumed to satisfy the strong slope conjecture, there are essential surfaces $S_1$ and $S_2$ in the exterior of $K$ such that 
\begin{enumerate}
\item the boundary slope of $S_1$ is $pq$ and $\chi(S_1)=0;$ and
\item  the boundary slope of $S_2$ is zero and $\displaystyle{{\chi(S_2)\over {{\abs{\partial S_2}} }}}= pq-p+q.$
\end{enumerate}

Since $\chi(S_1)=0$, by passing to the orientable double cover of $S_1$ if necessary, we conclude that the exterior of $K$ contains an essential annulus of boundary slope $pq$. It follows that
$K$ is a cable $J_{p_1,q_1}$ of a knot $J$,  where $p_1q_1=pq<0$. This follows, for example, by \cite[Theorem 4.18]{Budney}.
Note that the cabling annulus of $K=J_{p_1,q_1}$ satisfies Conjecture \ref{SSC} for $d_-[J_{K}(n)]$.
 
Next we need to gain a better understanding of the surface $S_2$:  First note that since $p$ and $q$ are relatively prime  integers  we have $\chi(S_2)\ne 0$; thus $S_2$ is not an annulus or a M\"obius band.

Recall that the exterior of $K$ is the union of a cable space $C_{p_1, q_1}$ and the exterior of $J$, say $M_J=\overline{S^3\setminus n(J)}$.  
Moreover $C_{p_1, q_1}$ is a  manifold with two torus boundary components, each of which is incompressible in $C_{p_1, q_1}$, that is a Seifert fibered manifold over an annulus with one singular fiber of order 
$q_1$. See \cite{Effie-Anh-slope}. By  classification of Seifert fibered manifold results (see, for example, \cite[Proposition 2.1 and Theorem 2.3]{Ha}), we can assume that $\abs{p_1}<\abs{q_1}$ and $q_1>0$.  Also since
replacing $(p_1, q_1)$ by  $(-p_1, -q_1)$ only results in knot orientation change that doesn't affect the colored Jones polynomial, we will assume that
$p_1<0<q_1$.
\vskip 0.05in 

{\it Claim:} $J$ is the trivial knot.

 \vskip 0.05in 

{\it Proof of Claim:}  The exterior  $M_K=\overline{S^3\setminus n(K)}$ is obtained by glueing together $M_J$ and $C_{p_1, q_1}$ 
along  a torus $T$. Note that if $T$ compresses in $M_K$ then $J$ is the trivial knot and we are done.

Suppose that $T$ is incompressible in $M_K$. Recall that we have an essential surface $S_2$ in $M_K$ such that
the boundary slope of $S_2$ is zero and $\displaystyle{{\chi(S_2)\over {{\abs{\partial S_2}} }}}= pq-p+q.$
We can isotope $S_2$ so that $S_2\cap T$ is a collection of essential, parallel curves in $T$.
As in the proof of \cite[Corollary 2.8]{Effie-Anh-slope} we have the following:

\begin{itemize}
\item The intersection $M_J\cap S_2$  is an essential surface  $S_J$ in $M_J$ such that each component of $\partial(M_J\cap S_2)$ has zero slope.

\item We have $\abs{\partial S_2}=\abs{\partial S_J}$.

\item We have

$${\chi(S_2)\over {{\abs{\partial S_2}} }}=q_1 {\chi(S_J)\over {{\abs{\partial S_J}} }}+ (q_1-1)p_1.$$
\end{itemize}
 
 Thus we have
 $$q_1 {\chi(S_J)\over {{\abs{\partial S_J}} }}+ (q_1-1)p_1=pq-p+q.$$
 
Since $p_1q_1=pq$  and $p<0<q$, we get
$$q_1 {\chi(S_J)\over {{\abs{\partial S_J} }}}-p_1=q-p>0,$$
which implies
$\displaystyle{{\chi(S_J)\over{{\abs{\partial S_J} }}}> p_1/q_1 > -1}$. The last  inequality gives $\chi(S_J)=2-2g(S_J)-{\abs{\partial S_J}} >-{\abs{\partial S_J}} $, where
$g(S_J)$ is the genus of $S_J$. Thus we must have $g(S_J)=0$
and $S_J$ is a collection of discs,  contradicting  our assumption that $T$ is incompressible. This finishes
the proof of Claim.
 \vskip 0.05in

Since $J$ is the unknot,  $K$ is the torus knot $T_{p_1, q_1}$ or $T_{-p_1, -q_1}$. However, since $p_1q_1=pq$ and $pq-p+q=p_1q_1-p_1+q_1,$
we get $p_1^2+q_1^2 =p^2+q^2$ and hence $(p_1+q_1)^2=(p+q)^2$. Now it follows easily 
 that $K$ is $T_{p, q}$ or $T_{-p, -q}$.
\qed

\vskip 0.06in

Theorem \ref{main}, and its proof, immediately give the following.
\begin{corollary} Suppose that $K$ satisfies Conjecture \ref{SSC}. If $d_+[J_{K}(n)]=d_+[J_{\rm{unknot}}(n)]$ or   $d_-[J_{K}(n)]=d_-[J_{{\rm{unknot}}}(n)]$, for all $n$,
then $K$ is isotopic to the unknot.
\end{corollary}

 \vskip 0.07in

\begin{figure}
  \includegraphics[scale=1.4]{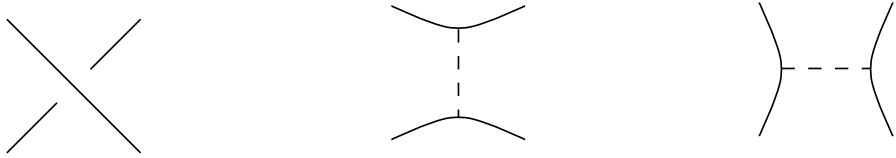}
  
  \caption{From left to right: A crossing, the $A$-resolution and  
 the $B$-resolution.}
    
  \label{resolutions}
\end{figure}

Next we discuss adequate knots: Let $D$ be a link diagram, and $x$ a crossing of $D$.  Associated to
$D$ and $x$ are two link diagrams,  called the \emph{$A$--resolution} and \emph{$B$--resolution} of
the crossing. See Figure \ref{resolutions}. Applying the $A$--resolution (resp. $B$--resolution)  to each crossing leads to  a collection of disjointly embedded circles $s_A(D)$
(resp.  $s_B(D)$).  The diagram $D$ is called  \emph{$A$--adequate}  (resp. \emph{$B$--adequate})
if for each crossing of $D$ the two arcs of $s_A(D)$
(resp.  $s_B(D)$) resulting from the resolution of the crossing lie on different circles. A knot diagram $D$
is \emph{adequate} if it is
both $A$-- and $B$--adequate.  Finally, a knot that admits an adequate diagram
is also called \emph{adequate}.

Starting with  $s_A(D)$
(resp.  $s_B(D)$) 
 we construct a surface $S_A(D)$ (resp.  $S_B(D)$)
as follows. Each circle of  $s_A(D)$
(resp.  $s_B(D)$)  bounds a disk on the projection sphere  $S^2\subset S^3$. This
collection of disks can be disjointly embedded in the ball below the
projection sphere. At each crossing of $D$, we connect the pair of
neighboring disks by a half-twisted band to construct a surface
whose boundary is $K$.  For details see \cite{FKP, FKP-guts}.
In \cite{ozawa} it is  shown that if $D$ is an adequate diagram of a knot $K$ then,
$S_A(D)$  and $S_B(D)$ are essential surfaces in the exterior of $K$. For a different proof of this fact see \cite{FKP-guts}.

\vskip 0.13in

  { \it Proof of Corollary \ref{adequate}.} Suppose that $K$ is an adequate knot such that  $d_+[J_{K}(n)]=d_+[J_{T_{p,q}}(n)]$ and  $d_-[J_{K}(n)]=d_-[J_{T_{p,q}}(n)]$, for all $n$.
  By Theorem \ref{main} $K$ is a torus knot. Thus we have a torus knot that is also adequate. Since $K$ is adequate,
there is an adequate diagram $D$.

It is known that  the number of negative crossings  $c_-(D)$  of an $A$--adequate knot diagram is a knot invariant.
Similarly, the number of positive crossings $c_+(D)$ of  a  $B$-adequate knot diagram is a knot invariant.
In fact, if $K$ is adequate, then the crossing number of $K$ is realized by the adequate diagram; that is  we have $c(K)=c(D)=c_-(D)+c_+(D)$ \cite{Li}.
Let $v_A(D)$ and  $v_B(D)$ be the numbers of  circles in $s_A(D)$ and $s_B(D)$, respectively.  The boundary slope of $S_A$ is $ -2c_- (D)$ and
$\chi(S_A)= v_A(D)-c(D)$.  The boundary slope of $S_B$ is $ 2c_+ (D)$ and
$\chi(S_B)=v_B(D)-c(D)$. 
By \cite{FKP}, the surfaces $S_A=S_A(D)$ and $S_B=S_B(D)$ satisfy the strong slope conjecture for $K$.

 That is 
we have
$$4\, d_-[J_{K}(n)] =  -2c_- (D) n^2 + 2(c(D) -v_A(D)) n  +2 v_A(D) -2 c_+(D),$$

and

$$
4 \, d_+[J_{K}(n)] = 2c_+ (D)n^2 + 2(v_B(D) - c(D)) n +2 c_-(D)-  2v_B(D).
$$

On the other hand,  $K$ is a torus knot. Without loss of generality we may assume  $K=T_{p,q}$, for $p<0<q$.
Thus $ -2c_-(D)=pq$ and $v_A(D) - c(D)=0$. This means that the surface $S_A(D)$ is a M\"obius strip. A spine of this surface is the graph
$\GA$ that is obtained as follows: the vertices are the circles of $v_A(D)$ and the edges are in one to one correspondence with the crossings
of $D$ in the way indicated by dashed lines in Figure \ref{resolutions}. The graph $\GA$ is connected and it  must retract to a circle. Thus it consists of circle $C$ with a number of trees attached to it. Consider a tree portion $T$ of  $\GA$ such that if $T$ is removed the Euler characteristic of the remaining graph is unchanged.
If $T$ contains edges then the diagram $D$ must contain nugatory crossings: For $T$  would define a connect summand of $S_A$ that is a disc; then each crossing corresponding to an edge of $T$ is a separating arc on the disc and thus the corresponding crossing of $D$ is nugatory (compare, proof of \cite[Corollary 3.21]{FKP-guts}). However since $D$ is adequate it cannot contain any nugatory crossings. Thus $T$ cannot contain any edges and $\GA$ consists of the circle $C$.
Now  each vertex of the graph $\GA$ 
is connected to exactly two edges. It follows that $D$ is the standard diagram of a $(-2,q)$-torus knot.
\qed

\smallskip

\begin{remark}{\rm   The proof of Corollary \ref{adequate} also shows  that the only adequate torus knots are those represented as closed 2-braids.  It may worth noting that this fact can be
deduced only using the quadratic terms of   $d_{\pm}[J_{T_{p,q}}(n)]$: The Jones slopes of   $T_{p,q}$ are $s^*=pq$ and $s=0$. If $T_{p,q}$ is adequate then $\abs{s-s^*}$ should be equal to twice the crossing number of $T_{p,q}$ (compare \cite[Theorem 1.1]{Kaind}). Thus we must have $\abs{pq}= \min \{ \abs{p(q-1)}, \abs{(p-1)q}\}$
which can happen only if $\abs{p}=2$  or $\abs{q}=2$.}
\end{remark}


\bibliographystyle{plain} \bibliography{biblio}

\begin{thebibliography}{10}

\bibitem{TakataMB}
Kenneth~L. Baker, Kimihiko Motegi, and Toshie Takata.
\newblock The strong slope conjecture for graph knots.
\newblock arXiv:1809.01039.

\bibitem{Takata}
Kenneth~L. Baker, Kimihiko Motegi, and Toshie Takata.
\newblock The strong slope conjecture for twisted generalized {W}hitehead
  doubles.
\newblock arXiv:1811.11673.

\bibitem{Budney}
Ryan Budney.
\newblock J{SJ}-decompositions of knot and link complements in {$S^3$}.
\newblock {\em Enseign. Math. (2)}, 52(3-4):319--359, 2006.

\bibitem{burde-zieschang:knots}
Gerhard Burde and Heiner Zieschang.
\newblock {\em Knots}, volume~5 of {\em de Gruyter Studies in Mathematics}.
\newblock Walter de Gruyter \& Co., Berlin, second edition, 2003.

\bibitem{FKP}
David Futer, Efstratia Kalfagianni, and Jessica~S. Purcell.
\newblock Slopes and colored {J}ones polynomials of adequate knots.
\newblock {\em Proc. Amer. Math. Soc.}, 139:1889--1896, 2011.

\bibitem{FKP-guts}
David Futer, Efstratia Kalfagianni, and Jessica~S. Purcell.
\newblock {\em Guts of surfaces and the colored {J}ones polynomial}, volume
  2069 of {\em Lecture Notes in Mathematics}.
\newblock Springer, Heidelberg, 2013.

\bibitem{ga-quasi}
Stavros Garoufalidis.
\newblock The degree of a {$q$}-holonomic sequence is a quadratic
  quasi-polynomial.
\newblock {\em Electron. J. Combin.}, 18(2):Paper 4, 23, 2011.

\bibitem{ga-slope}
Stavros Garoufalidis.
\newblock The {J}ones slopes of a knot.
\newblock {\em Quantum Topol.}, 2(1):43--69, 2011.

\bibitem{LeeGarVeen}
Stavros Garoufalidis, Christine Ruey~Shan Lee, and Roland van~der Veen.
\newblock The slope conjecture for {M}ontesinos knots.
\newblock arXiv:1807.00957.

\bibitem{Ha}
Allen Hatcher.
\newblock Notes on basic 3-manifold topology.
\newblock {\tt http://www.math.cornell.edu/\allowbreak \~{
  }hatcher/3M/\allowbreak 3Mdownloads.html}.

\bibitem{hatcher}
Allen~E. Hatcher.
\newblock On the boundary curves of incompressible surfaces.
\newblock {\em Pacific J. Math.}, 99(2):373--377, 1982.

\bibitem{Kaind}
Efstratia Kalfagianni.
\newblock A {J}ones slopes characterization of adequate knots.
\newblock {\em Indiana Univ. Math. J.}, 67(1):205--219, 2018.

\bibitem{Effie-Anh-slope}
Efstratia Kalfagianni and Anh~T. Tran.
\newblock Knot cabling and degrees of colored {J}ones polynomials.
\newblock {\em New York Journal of Mathematics}, Volume 21:905--941, 2015.

\bibitem{LeeVeen}
Christine Ruey~Shan Lee and Roland van~der Veen.
\newblock Slopes for pretzel knots.
\newblock {\em New York J. Math.}, 22:1339--1364, 2016.

\bibitem{Li}
W.~B.~Raymond Lickorish.
\newblock {\em An introduction to knot theory}, volume 175 of {\em Graduate
  Texts in Mathematics}.
\newblock Springer-Verlag, New York, 1997.

\bibitem{MoTa}
Kimihiko Motegi and Toshie Takata.
\newblock The slope conjecture for graph knots.
\newblock {\em Math. Proc. Cambridge Philos. Soc.}, 162(3):383--392, 2017.

\bibitem{ozawa}
Makoto Ozawa.
\newblock Essential state surfaces for knots and links.
\newblock {\em J. Aust. Math. Soc.}, 91(3):391--404, 2011.

\end{thebibliography}
\end{document}